\documentclass[10pt, reqno]{amsart}

\usepackage{hyperref} 
\usepackage{amssymb, mathrsfs,bbold}
\usepackage{amsmath, amsthm}
\usepackage{enumerate}
\usepackage{dsfont}
\usepackage{color}
\usepackage[a4paper]{geometry}
\usepackage[utf8]{inputenc}   
\usepackage{stmaryrd}
\usepackage{titlesec, wasysym}
\usepackage{appendix}
\usepackage{todonotes, fancyhdr}
\usepackage{chngcntr}
\usepackage{cancel}
\usepackage{mathtools}

\usepackage{tikz-cd}
\usetikzlibrary{calc}

\usepackage{setspace}
\renewcommand{\baselinestretch}{0.99}

\usepackage[all]{xy}
\numberwithin{subsection}{section}
\numberwithin{subsubsection}{subsection}
\numberwithin{equation}{section} 

\newenvironment{Dem}[1][\unskip]{%
    \begin{list}{\hspace{1.15cm}{\sf \textbf{{\small Proof #1 --}}}}{%
        \setlength{\topsep}{0pt}%
        \setlength{\leftmargin}{0pt}%
        \setlength{\rightmargin}{0pt}%
        \setlength{\listparindent}{0pt}%
        \setlength{\itemindent}{0pt}%
        \setlength{\parsep}{0pt}%
        \addtolength{\leftmargin}{20pt}%
        \addtolength{\rightmargin}{0pt}%
    } \item }{\hfill $\rhd$\end{list}\smallskip}

\newenvironment{Dem*}[1][\unskip]{%
    \begin{list}{\hspace{0cm}{\sf \textbf{{\small Proof #1 --}}}}{%
        \setlength{\topsep}{0pt}%
        \setlength{\leftmargin}{0pt}%
        \setlength{\rightmargin}{0pt}%
        \setlength{\listparindent}{0pt}%
        \setlength{\itemindent}{0pt}%
        \setlength{\parsep}{0pt}%
        \addtolength{\leftmargin}{20pt}%
        \addtolength{\rightmargin}{0pt}%
    } \item }{\hfill $\rhd$\end{list}\smallskip}

\newenvironment{Rem}[1][\unskip]{%
    \begin{list}{\hspace{1.15cm}{\textsf{\textbf{{\small \textsl{Remark}} --}}}}{%
        \setlength{\topsep}{0pt}%
        \setlength{\leftmargin}{0pt}%
        \setlength{\rightmargin}{0pt}%
        \setlength{\listparindent}{0pt}%
        \setlength{\itemindent}{0pt}%
        \setlength{\parsep}{0pt}%
        \addtolength{\leftmargin}{20pt}%
        \addtolength{\rightmargin}{0pt}%
    } \item }{\hfill \end{list}\smallskip}

\renewcommand\thesection       {\arabic{section}}
\renewcommand\thesubsection    {\thesection.\arabic{subsection}}
\renewcommand\thesubsubsection    {\thesection.\arabic{subsection}.\arabic{subsubsection}}

\titleformat{\section}[block]
{\filcenter\normalfont\sffamily\bfseries\Large}
{{\hspace{-0.7cm}}\thesection \hspace{0.2em} --\vspace{0.3cm}}{0.5em}{}

\titleformat{\subsection}[block]
{\filcenter\normalfont\sffamily\bfseries\large}  						  
{\hspace{-0.7cm}\thesubsection \hspace{0.5em} \vspace{0.3cm}}{.5em}{}  
\titlespacing{\subsection}{-0pc}{1.5ex plus .1ex minus .2ex}{0pc}

\titleformat{\subsubsection}[block]
{\normalfont\sffamily\bfseries}					  
{\thesubsubsection \vspace{0.3cm}}{.5em}{}  
\titlespacing{\subsection}{-0pc}{1.5ex plus .1ex minus .2ex}{0pc}

\newtheoremstyle{mystyle}
{3pt}               
{3pt}               
{\it }                      
{}                      
{\sffamily\bfseries}             
{}                      
{0.5em}                 
{#1 #2{\hspace{0.2cm}--\hspace{-0.2cm}}  }   

\theoremstyle{mystyle}

\newtheorem{thm}{Theorem}
\newtheorem*{thm*}{Theorem}

\newtheorem{cor}[thm]{\hspace{-0.15cm}  {Corollary} }
\newtheorem{lem}[thm]{\hspace{-0.14cm}  {Lemma} }
\newtheorem{prop}[thm]{\hspace{-0.13cm} {Proposition}}
\newtheorem{defn}[thm]{ \hspace{-0.31cm} {Definition}}

\newtheoremstyle{mystyle2}
{3pt}               
{3pt}               
{\it }                      
{}                      
{\sffamily\bfseries}             
{}                      
{0.5em}                 
{\llap{#2 }#1{\hspace{0.2cm}--}}

\theoremstyle{mystyle2}

\newtheorem*{definition*}{Definition}
\newtheorem*{theorem*}{Theorem}
\newtheorem*{Remark*}{Remark}
\newtheorem*{lem*} {Lemma}
\newtheorem*{defn*} {Definition}
\newtheorem*{prop*} {Proposition}
\newtheorem*{cor*} {Corollary}

\newcommand{\IR}{\mathbb{R}}

\newcommand{\IT}{\mathbb{T}}

\newcommand{\mcS}{\mathcal S}

\newcommand{\CH}{\mathcal H}

\newcommand{\ssk}{\smallskip}

\newcommand\bbE{\mathbb{E}}

\newcommand\bbR{\mathbb{R}}
\newcommand{\bbT}{\mathbb{T}}

\newcommand{\mcB}{\mathcal{B}}

\newcommand{\scrE}{\ensuremath{\mathscr{E}}}

\makeatletter
\newcommand*{\defeq}{\mathrel{\rlap{%
                     \raisebox{0.3ex}{$\m@th\cdot$}}%
                     \raisebox{-0.3ex}{$\m@th\cdot$}}%
                     =}
\makeatother

\begin{document}

\begin{center}
{\LARGE\sffamily{Variational methods for some singular stochastic elliptic PDEs  \vspace{0.5cm}}}
\end{center}

\begin{center}
{\sf I. BAILLEUL}\footnote{I.B. acknowledges support from the CNRS/PIMS and the ANR-16-CE40-0020-01 grant.} \& {\sf H. EULRY} \& {\sf T. ROBERT}
\end{center}

\vspace{1cm}

\begin{center}
\begin{minipage}{0.8\textwidth}
\renewcommand\baselinestretch{0.7} \scriptsize \textbf{\textsf{\noindent Abstract.}} We use some tools from nonlinear analysis to study two examples of stochastic elliptic PDEs, involving a singular operator, that cannot be solved by the contraction principle or the Schauder fixed point theorem. Let $\xi$ stand for a spatial white noise on a closed Riemannian surface $\mcS$. We prove the existence of a solution to the equation
$$
(-\Delta + a)u = f(u) + \xi u
$$
with a potential $a\in L^p(\mcS)$ and $p>1$, and $f$ subject to growth conditions. Under an additional parity condition on $f$ -- met for instance when $f(u) = u\vert u\vert^\ell$, with $\ell$ an even integer, we further prove that this equation has infinitely many solutions, in stark contrast with all the well-posedness results that have been proved so far for such singular stochastic PDEs under a small parameter assumption. This kind of results is obtained by seeing the equation as characterizing the critical points of an energy functional based on the Anderson operator $H=\Delta + \xi$ and by resorting to variants of the mountain pass theorem. There are however some interesting equations that cannot be characterized as the critical points of an energy functional. Such is the case of the singular Choquard-Pekar equation on $\mcS=\bbT^2$
$$
(-\Delta + a)u = \big(w\star f(u)\big)g(u) + \xi u
$$
One can use Ghoussoub's machinery of self-dual functionals to prove the existence of a solution to that equation as the minimum of a self-dual strongly coercive functional under proper assumptions on the coefficients $a,w,f$ and $g$.
\end{minipage}
\end{center}

\bigskip

\begin{center}
\begin{minipage}{0.8\textwidth}
\renewcommand\baselinestretch{0.7} \scriptsize \textbf{\textsf{\noindent R\'esum\'e.}} Nous utilisons des outils d'analyse non lin\'eaire pour \'etudier deux exemples d'\'equations aux d\'eriv\'ees partielles stochastiques (EDPS) elliptiques mettant en jeu un op\'erateur singulier, et qu'on ne peut r\'esoudre \`a l'aide d'une m\'ethode de point fixe classique, contraction ou Schauder. Soit $\xi$ un bruit blanc spatial d\'efini sur une surface riemannienne $\mcS$. Nous d\'emontrons l'existence d'une solution \`a l'\'equation
$$
(-\Delta + a)u = f(u) + \xi u
$$
o\`u $a$ est un potentiel $L^p(\mcS), p>1$, et la fonction $f$ satisfait de conditions de croissance. Sous une hypoth\`ese additionnelle de parit\'e, satisfaite par exemple lorsque $f(u) = u\vert u\vert^\ell$ avec $\ell$ un entier pair, nous d\'emontrons que l'\'equation admet une infinit\'e de solutions. Ce r\'esultat contraste fortement avec tous les r\'esultat d'existence d'une unique solution d\'emontr\'es jusqu'\`a pr\'esent pour toutes les EDPS singuli\`eres, sous des hypoth\`eses de petits param\`etres. Nos r\'esultats sont obtenus en caract\'erisant une solution comme point critique d'une fonctionnelle d'\'energie construite \`a partir de l'op\'erateur d'Anderson $H=\Delta + \xi$, et en faisant appel \`a des variantes du th\'eor\`eme du col. Cependant un certains nombre d'\'equations ne peuvent pas se formuler comme caract\'erisations de points critiques de fonctionnelles d'\'energies. C'est le cas de la version singuli\`ere de l'\'equation de Choquard-Pekar sur le tore $2$-dimensionel
$$
(-\Delta + a)u = \big(w\star f(u)\big)g(u) + \xi u,
$$
o\`u $\star$ d\'esigne l'op\'eration de convolution. On peut utiliser la machinerie des fonctionnelles auto-duales de Ghoussoub pour obtenir l'existence d'une solution \`a cette \'equation, sous la forme d'un point atteignant le minimum d'une fonctionnelle auto-duale fortement coercive sous certaines hypoth\`eses sur les coefficients $a,w,f,g$ de l'\'equation.
\end{minipage}
\end{center}

\vspace{0.6cm}

\section{Introduction}
\label{SectionIntro}

Let $(\mcS,g)$ stand for a closed (compact, connected, boundaryless) two dimensional Riemannian manifold. A spatial white noise $\xi$ on $\mcS$ is a random distribution with centered Gaussian law with covariance 
$$
\bbE\big[\xi(f_1)\,\xi(f_2)\big] = \int_\mcS f_1(x) f_2(x)dx
$$
for all smooth real-valued functions $f_1,f_2$ on $\mcS$, with $dx$ standing for the Riemannian volume measure. This random distribution takes almost surely its values in the Besov-H\"older space $B^{\alpha-2}_{\infty,\infty}(\mcS)$, for any $\alpha<1$ -- think of $\alpha-2$ as $(-1)^-$. The {\it Anderson operator} is formally defined as 
$$
Hu \defeq \Delta u+ \xi u,
$$
with $\xi$ seen here as a multiplication operator by $\xi$. The low regularity of $\xi$ causes problems to define $H$ as an unbounded operator on $L^2(\mcS)$ and is precisely what makes the equation singular. For the {\it product} $\xi u$ of $\xi$ by a function $u$ to make sense the function $u$ needs to have regularity $\beta$ with $(\alpha-2)+\beta>0$, that is $\beta=1^+$. The distribution $\Delta u$ will then have regularity $(-1)^+$, from which we should not expect any compensation with the $(-1)^-$ regularity of $\xi u$ to get an element $Hu$ of $L^2(\mcS)$ in the end. This state of affair can be disentangled using the tools of paracontrolled calculus or regularity structures to define $H$ as an unbounded operator on $L^2(\mcS)$ with a random domain $\mathfrak{D}(H)$. These tools have been developed for the study of singular stochastic partial differential equations after the pioneering works of Gubinelli, Imkeller \& Perkowski \cite{GIP} and M. Hairer \cite{Hai14}. The construction of the operator over a two dimensional torus was first performed by Allez \& Chouk in \cite{AllezChouk} using paracontrolled calculus. Their approach was subsequently simplified by Gubinelli, Ugurcan \& Zacchuber in \cite{GUZ}. Labb\'e gave the first construction of the Anderson operator over a three dimensional torus in \cite{Labbe} using the tools of regularity structures. Mouzard further simplified the approach of \cite{GUZ} and constructed the operator over an arbitrary closed two dimensional Riemannian manifold. A deep study of the Anderson operator and some associated objects was done recently by Bailleul, Dang \& Mouzard in \cite{BDM}. In any case one is able to define $H$ as a closed symmetric unbounded operator with random domain and compact resolvent. As such it has a nice spectral theory.

\ssk

We study in this work two classes of singular stochastic elliptic equations with multiplicative spatial white noise and prove existence results for them in settings where one cannot use a fixed point formulation of the equations. We are even able to exhibit a class of equations that have infinitely many solutions. This comes in stark contrast with all the well-posedness results proved in the literature on such singular stochastic partial differential equations under a small parameter assumption. This typically takes the form of existence (and uniqueness) for small times in the case of parabolic equations, e.g. \cite[Corollary 9.3]{Hai14} or \cite[Theorem 5.4]{GIP}, and small noise or strict convexity of the nonlinearity, as in \cite[Theorem 1.1]{OttoWeber} or \cite[Theorem 3]{ADG} for elliptic equations. We prove our results using a setting where solutions are understood in a weak sense and by resorting to variants of the mountain pass theorem. The use of topological methods to get critical points of $C^1$ functionals provides a very efficient and robust approach. There are however interesting equations that cannot be written as the Euler-Lagrange equation of some functional. The use of Ghoussoub's notion of self-dual functional provides a setting to characterize solutions of a number of equations as minimizers of a large class of functionals. Tools from convex analysis are required to set the scene. Note that the recent work \cite{DuanZhang} was the first to implement the variational method for the construction of solutions to elliptic equations associated with the Anderson operator, their approach relies on a direct method of calculus of variations and aims at proving some regularity results on said solutions.

\medskip

Section \ref{SectionBasicsAnderson} recalls and proves all we need to know about the Anderson operator and its perturbations by $L^p$ potentials. Section \ref{SectionWeakSolutions} is dedicated to the study of the equation
\begin{equation} \label{EqModelEq}
-Hu = au + f(\cdot,u),
\end{equation}
with potentials $a\in L^p(\mcS)$ for some $p>1$, with Theorem \ref{link_PDE} and Theorem \ref{fountain_PDE} as our main results. The statement gives mild conditions under which equation \eqref{EqModelEq} has at least one weak solution.  The second statement shows that an additional parity condition on the nonlinearity $f$ entails the existence of infinitely many weak solutions. Section \ref{SectionNonVariational} is dedicated to the study of the non-variational singular Choquard-Pekar equation on $\mcS=\bbT^2$
$$
-Hu = au + \big(w\star \vert u\vert^p\big) \vert u\vert^{q-2}u,
$$
for $p\neq q$. We obtain an existence result in Theorem \ref{ThmChoquardPekar} for some appropriate $w$, $p$ and $q$.

\bigskip

\noindent \textbf{\textsf{Notation --}} {\it All integrals will be with respect to the Riemannian volume measure. We will generically write them either as $\int_\mcS f$ or $\int_\mcS f(x)\,dx$.   }

\bigskip

\section{Basics on the Anderson operator}
\label{SectionBasicsAnderson}

We recall in Subsection \ref{SubsectionBasics} a number of results about the Anderson operator and prove in Subsection \ref{SubsectionKato} that the quadratic form associated with the Schr\"odinger Anderson operator $H+a$ has a nice spectral theory.

\medskip

\subsection{Basic results}
\label{SubsectionBasics}

We will not need in the present work any of the technical details associated with the use of paracontrolled calculus or regularity structures. We only mention from the works \cite{Mouzard, BDM} the following facts that we will freely use below, these are the only properties of the operator we need to conduct our analysis. Recall $\alpha-2<-1$ stands for the almost sure regularity of white noise.
\begin{itemize}
   \item[--] The Anderson operator $H$ can be defined as a closed symmetric unbounded operator on $L^2(\mcS)$ with random domain $\mathfrak{D}(H)$ and compact resolvent. As such it has a nice spectral theory and $H : \mathfrak{D}(H)\rightarrow L^2(\mcS)$ is almost surely invertible. (See e.g. Section 2 of \cite{Mouzard} or Section 3 of \cite{BDM}.)  \vspace{0.1cm}

   \item[--] There exists a random constant $c$ such that the quadratic form associated with the operator $-H+c$ is positive definite. The closure of the domain $\mathfrak{D}(H)$ with respect to the norm
   $$
   \Vert u\Vert_{\mathscr{E}} \defeq \sqrt{\big\langle(-H+c)u,u\big\rangle_{L^2}}
   $$
   defines a Hilbert space $\mathscr{E}$. That space is included and dense in any Sobolev space $H^{\beta}(\mcS)$, for $0\leq \beta<\alpha$, with compact inclusions. Set
   $$
   -H_c \defeq -H + c.
   $$

   \item[--] The operator $e^{tH_c}$ has a positive kernel $p_t(x,y)$ and there exists positive (random) constants $a_1,a_2$ such that one has
   \begin{equation} \label{EqGaussianBoundsHeatKernel}
   \frac{1}{a_1t}\,\exp\bigg(\hspace{-0.1cm}-a_2\frac{d(x,y)^2}{t}\bigg) \leq p_t(x,y) \leq \frac{a_1}{t}\,\exp\bigg(\hspace{-0.1cm}-\frac{d(x,y)^2}{a_2t}\bigg),
   \end{equation}
   uniformly in $x,y\in\mcS$ and $t\in(0,1]$, where $d(x,y)$ stands for the geodesic distance on $\mcS$ associated with the metric $g$. (See Proposition 25 in Section 4.3 of \cite{BDM}.)   \vspace{0.1cm}

   \item[--] There exists a positive (random) constant $\varepsilon$ such that 
   \begin{equation}\label{epsilon}
   e^{tH_c}{\bf 1} \leq e^{-t\varepsilon},
   \end{equation}
    for all $t>0$, and the Green function $G(x,y)$ of $H_c$ is finite outside the diagonal and satisfies the estimate
   \begin{equation} \label{EqEquivalenceGreenH}
  \big| \ln d(x,y)\big|\lesssim G(x,y) \lesssim \big\vert\hspace{-0.03cm}\ln d(x,y)\big\vert.
   \end{equation}
   (See\footnote{Strictly speaking, only an upper bound is stated in \cite[Lemma 36]{BDM}, but the equation (4.21) in \cite{BDM} shows that $(H+c)^{-1}=(\Delta+1)^{-1}+$regularizing operator, the latter having a continuous kernel on $\mcS\times\mcS$. Since the kernel of $(\Delta+1)^{-1}$ is exactly $-\frac1{2\pi}\ln d(x,y)+$ a continuous function on $\mcS\times\mcS$, both the upper and the lower bounds follow.} Proposition 25 in Section 4.3 and Lemma 36 in Section 5 of \cite{BDM}.)
\end{itemize}

\medskip

We do not record in the heat kernel $p_t$ or the Green function $G$ the dependence of these functions on the constant $c$ as the latter will be fixed throughout. It follows from the second item and the Sobolev embedding that we have a compact inclusion of $\mathscr{E}$ into $L^q(\mcS)$, for all $1<q<\frac{2}{1-\alpha}$. Any bounded sequence in $\mathscr{E}$ has thus a subsequence that converges weakly in $\mathscr{E}$ and strongly in $L^q(\mcS)$, for a given $1<q<\frac{2}{1-\alpha}$. (We will use that fact a few times.) Do not be mislead by the comparison of the Green function of $H$ with the Green function of $\Delta$ in the fourth item. While we have the small distance bound \eqref{EqEquivalenceGreenH} between the two functions the integral operator on functions associated with $G$ does not have the regularizing properties that the operator $\Delta^{-1}$ have: there is {\it no} elliptic regularity for the operator $H^{-1}$. This fact is related to the singular character of the Anderson operator and the low regularity of white noise.

\medskip

It is already possible from these facts to say something about the solvability of the semilinear stationary Schr\"odinger Anderson equation
\begin{equation} \label{EqSchrodingerAnderson}
-Hu = au + f(\cdot,u)
\end{equation}
when the right hand side is a priori in $L^2(\mcS)$, using the (almost sure) invertibility of $H$ and the compact embedding of its domain in $L^2(\mcS)$.

\medskip

\begin{prop} \label{PropSchauder}
Assume that $a\in L^\infty(\mcS)$ and that one can associate to $f\in C^0(\mcS\times\IR)$ a function $h\in L^2(\mcS)$ such that $\left|f(\cdot,z)\right|\leq h(\cdot)$, uniformly in $z\in\IR$. Then equation \eqref{EqSchrodingerAnderson} has a solution if $\|a\|_{L^\infty}$ is small enough.
\end{prop}

\medskip

\begin{Dem}
The continuity of the operator $H^{-1} : L^2(\mcS)\rightarrow L^2(\mcS)$ and the estimate
$$
\big\Vert au+f(\cdot,u)\big\Vert_{L^2} \leq \Vert a\Vert_\infty \Vert u\Vert_{L^2} + \Vert h\Vert_{L^2}
$$
tell us that a ball of $L^2(\mcS)$ of large enough radius is sent by the map $u\mapsto H^{-1}\big(au+f(\cdot,u)\big)$ into itself. As $H^{-1}$ actually takes values in the compact subset $\mathfrak{D}(H)$ of $L^2(\mcS)$ the conclusion comes from Schauder fixed point theorem.
\end{Dem}

\medskip

Alternatively, for $a\in L^2(\mcS)$ one can use the Cameron-Martin theorem to say that the operator has a law that is equivalent to the law of $H$. (One could even use the much refined form of Cameron-Martin theorem proved by Kusuoka for random potentials $a$, under appropriate assumptions -- see e.g. Theorem 3.5.4 in \cite{UstunelZakai}.) So the almost sure existence of a solution to equation \eqref{EqSchrodingerAnderson} is equivalent in that case to the almost sure existence of a solution to equation
$$
Hu = f(\cdot,u).
$$
One can use a Schauder fixed point strategy if $f$ satisfies for instance an estimate of the form
$$
\big\Vert f(\cdot,u)\big\Vert_{L^2} \lesssim 1+o(\Vert u\Vert_{L^2})
$$
when $\Vert u\Vert_{L^2}$ goes to $+\infty$. This is in particular the case when $\vert f(\cdot,z)\vert\lesssim 1+ \vert z\vert^{\ell-1}$, for $\ell\leq 2$. While the compactness/(fixed point) method is elementary to set up it requires in one form or another a small size or integrability assumption on $a$. The topological methods used in Section \ref{SectionWeakSolutions} will bypass that constraint and work without size conditions on $a$ for the much larger class of $L^p$ potentials, for any $p>1$. As a preliminary step to the developments of Section \ref{SectionWeakSolutions} we first study the Schr\"odinger Anderson operator
$$
u\mapsto (-H + a)u
$$
for itself and give conditions on the potential $a$ in the next section for its associated quadratic form to have a nice spectral theory. These conditions are met for a large class of potentials, including $a\in L^p(\mcS)$ when $p>1$.

\bigskip

\subsection{The Kato class and the Anderson operator}
\label{SubsectionKato}

The aim of this section is to prove the following diagonalisation result for the Schr\"odinger Anderson operator $-H+a$.

\medskip

\begin{thm} \label{ThmSpectralDecompositionAndersonPotential}
Pick $a\in L^p(\mcS)$ with $p>1$. There exists an orthonormal basis $(e_i)_{i\geq 0}$ of $L^2(\mcS)$ such that 
$$
\mathscr{E} = \overline{\bigoplus_{i\geq 0}\bbR e_i},
$$
with the closure in $\mathscr{E}$, and one has for all $i\geq 0$
$$
\big\langle e_i , (-H+a)e_i\big\rangle_{L^2} = \mu_i.
$$
\end{thm}

\medskip

Recall that a potential $a : \mcS\rightarrow\bbR$ is said to be in the {\it Kato class} if 
\begin{equation} \label{EqDefnKatoClass}
\lim_{r\to0^+} \, \sup_{x\in\mcS} \int_{d(x,\cdot)<r} \big\vert\hspace{-0.03cm}\ln d(x,y)\big\vert\,|a(y)|\,dy = 0.
\end{equation}
Note that as $\mcS$ is compact, potentials in $L^p(\mcS)$ with $p>1$ are in the Kato class and that Kato class potentials are integrable.

Given the equivalence \eqref{EqEquivalenceGreenH} for the Green function $G$ of the Anderson operator $H$ one can rewrite condition \eqref{EqDefnKatoClass} under the form
$$
\lim_{r\to0^+} \, \sup_{x\in\mcS} \int_{d(x,y)<r} G(x,y)|a(y)|\,dy = 0.
$$

The proof of Theorem \ref{ThmSpectralDecompositionAndersonPotential} follows the proof of a similar result for perturbations of the $\Delta$ operator by potentials in the Kato class. (See for instance Section 3.3 of the book \cite{BetzLorincziHiroshima} of Betz, Hiroshima \& Lorinczi.) We rewrite in Proposition \ref{Kato_class_carac} condition \eqref{EqDefnKatoClass} as a condition on the operator $-H_c+\lambda$, when the constant $\lambda$ goes to $\infty$, and deduce from it in Proposition \ref{PropSmallBoundA} that the quadratic form associated with $a$ is $(-H_c)$-form bounded with arbitrarily small relative bound. We first state and prove these two propositions before proving Theorem \ref{ThmSpectralDecompositionAndersonPotential}. An intermediate result is needed first.

\medskip

\begin{lem} \label{LemKato}
A function $a\in L^1(\mcS)$ is in the Kato class iff
\begin{equation} \label{EqConditionLemmaKato}
\sup_{x\in\mcS}\int_0^T\int_{\mcS}p_s(x,y)|a(y)|\,dyds \underset{T\rightarrow 0^+}{\longrightarrow} 0.
\end{equation}
\end{lem}

\medskip

\begin{Dem}
We first note from the Gaussian bounds \eqref{EqGaussianBoundsHeatKernel} that condition \eqref{EqConditionLemmaKato} is equivalent to the condition
\begin{equation} \label{EqConditionLemmaKatoReformulation}
\sup_{x\in\mcS}\int_0^T\int_{\mcS}s^{-1}e^{-d(x,y)^2/s}|a(y)\vert\,dyds \underset{T\rightarrow 0^+}{\longrightarrow} 0.
\end{equation}

\ssk

{\color{gray} $\bullet$} Let $a$ be a potential in the Kato class. For $0<T<1$, we split the integration over $\mcS$ in \eqref{EqConditionLemmaKatoReformulation} into $\{d(x,\cdot)<T^{1/4}\}\cup\{d(x,\cdot)\geq T^{1/4}\}$. By Fubini-Tonelli's theorem, a change of variables, and integration by parts, one has
\begin{equation*} \begin{split}
\int_0^T&\int_{d(x,\cdot)<T^{1/4}} s^{-1}e^{-d(x,y)^2/s}|a(y)\vert\,dyds = \int_{d(x,\cdot)<T^{1/4}}\int_{T^{-1}d(x,y)^2}^{+\infty} r^{-1}e^{-r}|a(y)\vert\,drdy   \\
&\lesssim -\int_{d(x,\cdot)<T^{1/4}} \ln\bigg(\frac{d(x,y)^2}{T}\bigg)|a(y)\vert\,dy + \int_{d(x,\cdot)<T^{1/4}}\int_{T^{-1}d(x,y)^2}^{+\infty} (\ln r)\,e^{-r}|a(y)\vert\,drdy   \\
&\lesssim \int_{d(x,\cdot)<T^{1/4}} \big\vert\hspace{-0.03cm}\ln d(x,y)\big\vert\,|a(y)|\,dy + \ln T + o_T(1),
\end{split} \end{equation*}
with a $o_T(1)$ that comes from the integrable character of $a$ and a negative contribution of $\ln T$ that can be skipped in an upper bound. 

As we also have 
\begin{equation*} \begin{split}
\int_0^T\int_{d(x,\cdot)\geq T^{1/4}} s^{-1}e^{-d(x,y)^2/s}|a|(y)\,dyds &=\int_{d(x,\cdot)\geq T^{1/4}} \int_{T^{-1}d(x,y)^2}^{+\infty} r^{-1}e^{-r}|a|(y)\,dydr   \\
&\leq \int_{d(x,\cdot)>T^{1/4}} \int_{T^{-1/2}}^{+\infty} r^{-1}e^{-r}|a|(y)\,drdy = o_T(1),
\end{split} \end{equation*}
from the fact that $a\in L^1(\mcS)$, we see that condition \eqref{EqConditionLemmaKatoReformulation} follows from condition \eqref{EqDefnKatoClass}.

\ssk

{\color{gray} $\bullet$} Write $A\asymp B$ when we have both $A\lesssim B$ and $B\lesssim A$. We have the estimate
\begin{equation*} \begin{split}
\int_0^T p_s(x,y)ds &\asymp \int_0^T s^{-1}e^{-d(x,y)^2/s}ds \asymp \int_{d(x,y)^2/T}^{+\infty} r^{-1}e^{-r}dr  \\
&\asymp -\ln\big(d(x,y)^2/T\big)\,e^{-d(x,y)^2/T} + \int_{d(x,y)^2/T}^{+\infty} (\ln r)\,e^{-r}\,dr,
\end{split} \end{equation*}

which holds for any $0<T<1$ and uniformly in $x,y\in\mcS$,
and thus the upper bound
$$
\big(-\ln d(x,y)\big){\bf 1}_{d(x,y)\leq T} \lesssim \int_0^T p_s(x,y)ds + {\bf 1}_{d(x,y)\leq T}.
$$
Multiplying by $|a|$, integrating on $\mcS$ and using again Fubini-Tonelli's theorem, we see on this inequality that condition \eqref{EqDefnKatoClass} follows from condition \eqref{EqConditionLemmaKatoReformulation}.
\end{Dem}

\medskip

\begin{prop} \label{Kato_class_carac}
A function $a\in L^1(\mcS)$ is in the Kato class iff $\left\|(-H_c+\lambda)^{-1}|a|\right\|_\infty \underset{\lambda\rightarrow+\infty}{\longrightarrow} \, 0.$
\end{prop}

\medskip

\begin{Dem}
While the operator $(-H_c+\lambda)^{-1}$ is first defined as an operator from $L^2(\mcS)$ into $\mathfrak{D}(H)$, for good $\lambda$'s, its spectral representation
$$
(-H_c+\lambda)^{-1}u = \int_0^{+\infty}\int_{\mcS} e^{-\lambda t}p_t(\cdot,x)u(x)\,dxdt
$$
allows to extend it naturally to the set of non-negative valued functions $u$, with $(-H_c+\lambda)^{-1}u$ taking values in $[0,+\infty]$. (Recall the heat kernel of $H_c$ is positive, so the above quantity is positive unless $u$ is null.) Take $T>0$ to be chosen later. Slicing the time integral and changing variables, we have
\begin{equation*} \begin{split}
\big((-H_c+\lambda)^{-1}\vert a\vert\big)(z) = \sum_{n\geq 0} e^{-T\lambda n} \int_0^Te^{-\lambda s}\int_\mcS p_{nT}(z,x)\big(e^{sH_c}\vert a\vert\big)(x)\,dxds.
\end{split}
\end{equation*}
Thus with $\varepsilon$ as in \eqref{epsilon}, we see from the fact that $e^{nTH_c}{\bf 1}\leq e^{-nT\varepsilon}$ and the spectral representation of $(-H_c+\lambda)^{-1}$ that one has the upper bound
\begin{equation*}\begin{split}
\big((-H_c+\lambda)^{-1}\vert a\vert\big)(z) \leq \frac{1}{1-e^{-(\lambda+\varepsilon)T}}\;\underset{x\in\mcS}{\sup}\;\int_0^T\big(e^{sH_c}\vert a\vert\big)(x)\,ds \\ \lesssim  \frac{e^{\lambda T}}{1-e^{-(\lambda+\varepsilon)T}}\,\big\Vert (-H_c+\lambda)^{-1}\vert a\vert\big\Vert_{\infty}.
\end{split}
\end{equation*}
Taking $T=1/(\lambda+\varepsilon)$ shows then that we have
$$
\big\Vert (-H_c+\lambda)^{-1}\vert a\vert\big\Vert_{\infty} \asymp \underset{x\in\mcS}{\sup}\;\int_0^T\big(e^{-sH_c}\vert a\vert\big)(x)\,ds.
$$
As 
$$
\int_0^Te^{-sH_c}|a|(x)ds = \int_0^T\int_{\mcS}p_s(y,x)|a|(y)\,dyds
$$
and $p_s(\cdot,\cdot)$ is a symmetric function of its two space arguments the quantity $\big\Vert (-H_c+\lambda)^{-1}\vert a\vert\big\Vert_{\infty}$ is equivalent to the quantity $\int_0^T\int_{\mcS}p_s(x,y)|a|(y)\,dyds$, so the conclusion follows from Lemma \ref{LemKato}.
\end{Dem}

\medskip

\begin{prop} \label{PropSmallBoundA}
Let $a$ be a potential in the Kato class. For any $\eta>0$ there exists a positive constant $m_\eta$ such that one has 
$$
\langle u,|a|u\rangle_{L^2} \leq \eta\|u\|_{\mathscr{E}}^2 + m_\eta\|u\|_{L^2}^2,
$$
for all $u\in\mathscr{E}$.
\end{prop}

\medskip

\begin{Dem}
We prove below that the operator $\vert a\vert^{1/2}(-H_c+\lambda)^{-1/2}$ is well defined as an operator from $L^2(\mcS)$ into itself, with operator norm of order $\big\Vert (-H_c+\lambda)^{-1}\vert a\vert\big\Vert_{\infty}^{1/2}$. The inequality of the statement then follows from the identity
\begin{eqnarray*}
\langle u,|a|u\rangle_{L^2}&=&\||a|^{1/2}u\|_{L^2}^2=\left\||a|^{1/2}(-H_c+\lambda)^{-1/2}(-H_c+\lambda)^{1/2}u\right\|_{L^2}^2\\
&\leq&\left\||a|^{1/2}(-H_c+\lambda)^{-1/2}\right\|_{L^2\to L^2}^2\left\|(-H_c+\lambda)^{1/2}u\right\|_{L^2}^2,
\end{eqnarray*}
valid for $u\in\mathscr{E}$, and Proposition \ref{Kato_class_carac}.

\ssk

Now note first that 
$$
\left\|(-H_c+\lambda)^{-1}|a|\right\|_{L^\infty\to L^\infty} = \left\|(-H_c+\lambda)^{-1}|a|\right\|_\infty.
$$
By duality $|a|(-H_c+\lambda)^{-1}$ defines a bounded operator from $L^1(\mcS)$ into itself, with operator norm
$$
\left\||a|(-H_c+\lambda)^{-1}\right\|_{L^1\to L^1} = \left\|(-H_c+\lambda)^{-1}|a|\right\|_{L^\infty\to L^\infty} = \left\|(-H_c+\lambda)^{-1}|a|\right\|_\infty.
$$
Stein's interpolation theorem can thus be applied to the holomorphic family of operators
$$
T(z) \defeq |a|^z(-H_c+E)^{-1}|a|^{1-z},
$$
and shows that $T(1/2)$ is a bounded operator from $L^2(\mcS)$ into itself with operator norm at most $\left\|(-H_c+\lambda)^{-1}|a|\right\|_\infty$. The conclusion follows then from the identity
$$
\left\||a|^{1/2}(-H_c+\lambda)^{-1/2}\right\|_{L^2\to L^2}^2=\left\||a|^{1/2}(-H_c+\lambda)^{-1}|a|^{1/2}\right\|_{L^2\to L^2}=\|T(1/2)\|_{L^2\to L^2}.
$$
\end{Dem}

\medskip

The statement of Theorem \ref{ThmSpectralDecompositionAndersonPotential} is then a direct consequence of classical results on perturbations of quadratic forms, as Proposition \ref{PropSmallBoundA} allows us to use Theorem X.17 and Theorem XIII.68 of Reed \& Simon's books \cite{ReedSimon2} and \cite{ReedSimon4}, respectively. We order the family of the real-valued (random) eigenvalues of the quadratic form $-H_c+a$
\begin{equation} \label{EqDefnIndexm}
\mu_0\leq\mu_1\leq\cdots\leq\mu_m\leq0<\mu_{m+1}\leq\cdots
\end{equation}
and denote by $\mu_{m+1}$ the smallest positive eigenvalue -- with the convention that $m=-1$ if $\mu_0>0$. We record here for later use the following elementary result. Set 
$$
\mathscr{E}_{>m} \defeq \overline{\bigoplus_{i\geq m+1} \bbR e_i},
$$
with closure in $\mathscr{E}$.

\medskip

\begin{lem} \label{LemDefnSmallDelta}
Let $a\in L^p(\mcS)$ for some $p>1$. Then the following quantity is positive
$$
\delta := \underset{\Vert v\Vert_\mathscr{E}=1}{\underset{v\in\mathscr{E}_{>m}}{\inf}}\; \Big(\Vert v\Vert_\mathscr{E}^2 + \int_\mcS av^2\Big)  > 0.
$$
\end{lem}

\medskip

\begin{Dem}
We use the fact that $\alpha<1$ can be chosen arbitrarily close to $1$ to pick it in such a way that $2p/(p-1)<2/(1-\alpha)$. Recall that the space $\mathscr{E}$ is compactly embedded in $H^\beta(\mcS)$ for any $0\le \beta <\alpha$. Take a minimizing sequence $u_n$ in $\mathscr{E}_{>m}$ with $\|u_n\|_{\mathscr{E}}=1$, such that $\|u_n\|_{\mathscr{E}}+\int_\mcS au_n^2 = 1+\int_\mcS au_n^2\to \delta$. Then, since the sequence $u_n$ is bounded in $\mathscr{E}$ and takes values in the closed subspace $\mathscr{E}_{>m}$ it has a subsequence that converges weakly to an element $u$ of $\mathscr{E}_{>m}$ and, together with Sobolev embedding, strongly to $u$ in $L^{2p/(p-1)}(\mcS)$.  The integrals $\int_\mcS au_n^2$ then converge to $\int_\mcS au^2$, and 
$$\delta = 1+\int_\mcS au^2 = \liminf_{n\to\infty}\|u_n\|_{\mathscr{E}}^2+\int_\mcS au^2 \ge\|u\|_{\mathscr{E}}^2+\int_\mcS au^2.$$
If $u=0$ we have $\delta=1$, otherwise since $u\in\mathscr{E}_{>m}$ we have
$$
\delta \ge \int_\mcS \Big(\big(\sqrt{-H_cu}\big)^2+au^2\Big) \geq \mu_{m+1}\Vert u\Vert_\mathscr{E}^2.
$$
\end{Dem}

\medskip

\begin{Rem}
\it Note that even if Theorem \ref{ThmSpectralDecompositionAndersonPotential} and Lemma \ref{LemDefnSmallDelta} would hold true for any Kato class potential $a$ as well, we shall only consider the case where $a\in L^p(\mcS)$ for some $p>1$ in the following. This is required for our energy functional $\Phi$ defined below to be $C^1$.
\end{Rem}

\bigskip

\section{Weak solutions to singular stochastic PDEs}
\label{SectionWeakSolutions}

Let a function $f : \mcS\times\bbR\rightarrow\bbR$ be given, with $f(x,\cdot)\in L^1_{\textrm{loc}}(\bbR)$ for each $x\in\mcS$ and 
$$
\big\vert f(x,z)\big\vert \lesssim 1+\vert z\vert^\ell,
$$
for some positive exponent $\ell$, uniformly in $x\in\mcS$. We associate to $f$ the function
\begin{equation} \label{EqDefnF}
F(x,z) \defeq \int_0^z f(x,r)\,dr,
\end{equation}
defined for all for $(x,z)\in\mcS\times\bbR$. Pick $a\in L^p(\mcS)$ with $p>1$ and set
$$
\Phi(u) \defeq \frac{1}{2}\Vert u\Vert_\mathscr{E}^2 + \int_\mcS\bigg(\frac{1}{2}\,a(x)u(x)^2 - F\big(x,u(x)\big)\hspace{-0.07cm}\bigg)dx.
$$

\medskip

\begin{lem}
The function $\Phi$ on $\mathscr{E}$ is well-defined and $C^1$, with Fr\'echet derivative 
$$
\Phi'(u)(v) = \langle u,v\rangle_{\mathscr{E}} + \int_{\mcS}\bigg(\hspace{-0.05cm}a(x)u(x)v(x) - f\big(x,u(x)\big)v(x)\hspace{-0.07cm}\bigg)dx
$$

\end{lem}

\medskip

\begin{Dem}
We use again the fact that $\alpha<1$ can be chosen arbitrarily close to $1$ to pick it in such a way that $2p/(p-1)<2/(1-\alpha)$ and $\ell+1<2/(1-\alpha)$. The continuous embedding of $\mathscr{E}$ into $L^{2p/(p-1)}(\mcS)$ then tells us that the integral $\int_\mcS au^2$ defines a $C^1$ function of $u\in\mathscr{E}$ with derivative $v\mapsto 2\int_\mcS auv$ at point $u$. Similar considerations give the Fr\'echet differentiability of $\int_\mcS F\big(x,u(x)\big)dx$ as a function of $u\in\mathscr{E}$ and the formula for its derivative.
\end{Dem}

\medskip

This result justifies the following definition.

\medskip

\begin{defn} \label{DefnWeakSolutions}
A \textbf{\textsf{weak solution of the equation}}
\begin{equation} \label{EqModelEquation}
-H_cu +au= f(\cdot,u)
\end{equation}
is a critical point of the map $\Phi$.
\end{defn}

\medskip
Note that as we are working in a Hilbert space framework one can identify the Fréchet derivative of $\Phi$ at point $u$ to its gradient in $\mathscr{E}$, still denoted by $\Phi'(u)$. In their recent work \cite{DuanZhang} Duan \& Zhang use a similar characterization of weak solutions to the same equation. However their approach only deals with constant potentials $a\equiv\mu$, for which the direct method in the calculus of variations applies since the energy functional is coercive. This also allows them to prove that their weak solution actually belong to $\mathfrak{D}(H)$ and even get Schauder estimates on said solution. However in our case the whole point is that even the existence of a minimizer is not straightforward, since the energy functional is no longer coercive when $a$ is not constant and non-positive.
We note also that Igant, Otto, Ried \& Tsatsoulis also used slightly earlier a variational characterization of solutions to a nonlocal singular equation in their work \cite{IORT}. In echo of Definition \ref{DefnWeakSolutions} we define a weak solution of the equation
\begin{equation*}
-Hu + au = f(\cdot,u)
\end{equation*}
as a critical point of the map on $\mathscr{E}$ constructed with $a-c$ in place of $a$ -- shifting the function $a$ by a constant keeps its integrability property as $\mcS$ has finite volume. Working with the positive definite operator $-H_c$ turns out to be more practical. Note that one cannot use any kind of bootstrap, or elliptic regularity result, to get that weak solutions of equation \eqref{EqModelEquation} are strong solutions of that equation, as this would require $a$ to be an element of $L^2(\mcS)$. Indeed, $u$ cannot be expected to have more than $\CH^{1^-}(\mcS)$ regularity as an element of the domain $\mathfrak{D}(H)$, thus for the product $au$ to be in $L^2(\mcS)$ the potential $a$ would need to have enough integrability, which we do not assume. Here our argument covers any potential $a\in L^p(\mcS)$ in the whole range $p>1$.

\bigskip

\subsection{The Mountain Pass strategy}
\label{SubsectionMoutainPass}

We use a well-known variant of the mountain pass theorem to guarantee the existence of critical points of $\Phi$ under appropriate assumptions on $f$. First recall the following definition.

\medskip

\begin{defn*}
Let $b\in\bbR$. The functional $\Phi$ is said to satisfy the \textbf{\textsf{Palais-Smale condition $(\textsf{PS})_b$}} if any sequence $(u_n)$ in $\mathscr{E}$ satisfying
\begin{equation} \label{EqPS}
\Phi(u_n)\underset{n\rightarrow+\infty}{\longrightarrow}b, \qquad \Phi'(u_n)\underset{n\rightarrow+\infty}{\longrightarrow}0,
\end{equation}
has a converging subsequence in $\mathscr{E}$.
\end{defn*}

\medskip

With this property, the mechanics of minimax principles is simple and can be illustrated on the following special case.

\ssk

Let $B$ stand for the closed unit ball of the $d$-dimensional Euclidean space. Let $\rho_0$ be a continuous map from the unit sphere $\partial B$ into $\mathscr{E}$. Let $\Gamma$ stand for the collection of all continuous maps from $B$ into $\mathscr{E}$ whose restriction to $\partial B$ is $\rho_0$. If
\begin{equation} \label{EqMountainPassCondition}
\underset{\vert z\vert=1}{\max} \, \Phi\big(\rho_0(z)\big) < b\defeq \underset{\rho\in\Gamma}{\inf}\, \Vert\Phi\circ\rho\Vert_{\infty} < \infty
\end{equation}
then one can associate to every $\theta>0$ and every $\rho\in\Gamma$ such that
$$
\Vert\Phi\circ\rho\Vert_{\infty} \leq b+\theta
$$
a point $u\in\mathscr{E}$ such that
\begin{equation*} \begin{split}
\begin{cases}
&\big\vert \Phi(u)-b\big\vert \leq 2\theta,   \\
&\textrm{dist}\big(u,\rho(B)\big) \leq 2,   \\
&\Vert\Phi'(u)\Vert \leq 8\theta.
\end{cases}
\end{split} \end{equation*}
Indeed if all points of the $2$-neighbourhood of $\rho(B)$ where $\big\vert \Phi(u)-b\big\vert \leq 2\theta$ satisfied $\Vert\Phi'(u)\Vert > 8\theta$ one could build an explicit deformation $\widetilde{\rho}$ of $\rho$ that would be in the family $\Gamma$ and would satisfy $\Vert\Phi\circ\widetilde{\rho}\Vert_{\infty} \leq b-\theta$, contradicting the definition of $b$. Such a deformation would be constructed from the flow of a pseudo-gradient vector field associated with $\Phi'$. See e.g. Lemma 2.2, Lemma 2.3 and Theorem 2.8 in Willem's book \cite{Willem} -- here we took $\delta=1$ in the notations of \cite{Willem}. So there exists a sequence of points $u_n\in\mathscr{E}$ satisfying
$$
\Phi(u_n)\underset{n\rightarrow+\infty}{\longrightarrow}b, \qquad \Phi'(u_n)\underset{n\rightarrow+\infty}{\longrightarrow}0.
$$
If $\Phi$ satisfies the $(\textsf{PS})_b$ condition, any limit point $u$ is thus a critical point of $\Phi$ where $\Phi(u)=b$.

Let $S_r\subset\mathscr{E}$ stand for the sphere of $\mathscr{E}$ of radius $r$. If the maps $\rho\in\Gamma$ are of the form $\overline{\rho}\circ\iota$, where $\iota$ sends homeomorphically $B$ into $\mathscr{E}$ and $\iota(B)\cap S_r\neq\emptyset$, with the maps $\overline{\rho}$ defined on $\iota(B)$, they satisfy $\rho(B)\cap S_r\neq\emptyset$ -- otherwise one could construct a continuous retraction from $B$ into $\partial B$. (See e.g. the proof of Theorem 2.12 in \cite{Willem}.) Condition \eqref{EqMountainPassCondition} thus holds true if 
$$
\underset{\vert z\vert=1}{\max} \, \Phi\big(\rho(z)\big) < \underset{S_r}{\inf}\, \Phi.
$$
The (slightly refined) form under which we will use that fact is given by \textbf{\textsf{Rabinowitz' linking theorem}}, which we formulate in our setting here; see e.g. \cite[Theorem 2.12]{Willem}. Set for all $k\geq 0$
$$
\mathscr{E}_{\leq k} \defeq \bigoplus_{i=0}^k \bbR e_i, \qquad \mathscr{E}_{>k} \defeq \overline{\bigoplus_{i\geq k+1} \bbR e_i},
$$
with closure in $\mathscr{E}$.

\medskip

\begin{thm} \label{linking}
Pick $0<r_1<r_2<\infty$ and ${\sf y}\in\mathscr{E}_{>k}$ with norm $r_1$. Set
\begin{equation*} \begin{split}
\mcB_{r_2} &\defeq \Big\{u=y+t{\sf y}\,,\ y\in\mathscr{E}_{\leq k}\,;\, t\geq0\ \text{such that}\ \|u\|\leq r_2\Big\},
\end{split} \end{equation*}
and let $\Gamma$ stand for the set of continuous maps from $\mcB_{r_2}$ into $\mathscr{E}$ whose restriction to $\partial\mcB_{r_2}$ is the identity map. Then 
$$
b \defeq \underset{\rho\in\Gamma}{\inf}\,\Vert \Phi\circ\rho\Vert_\infty
$$
is a critical value of $\Phi$ if $\Phi$ satisfies the Palais-Smale condition $(\textsf{\emph{PS}})_b$ and
\begin{equation} \label{EqRabinowitzCondition}
\underset{\partial\mcB_{r_2}}{\max}\,\Phi < \underset{S_{r_1}\cap \mathscr{E}_{>k}}{\inf}\,\Phi.
\end{equation}
\end{thm}

\medskip

We will use that result to prove existence of weak solutions of equation \eqref{EqModelEquation}. The following variation on Rabinowitz' linking theorem due to Bartsch will be used to prove that equation \eqref{EqModelEquation} actually have infinitely many solutions under an appropriate parity assumption on $f$.  Here again the statement is given in our setting, and we refer e.g. to \cite[Theorem 3.6]{Willem}.

\medskip

\begin{thm} \label{ThmBartsch} \textbf{\textsf{(Barstch's fountain Theorem)}}
Assume $f$ is odd with respect to its $z$ argument. If $\Phi$ satisfies the Palais-Smale condition $(\textsf{\emph{PS}})_b$ for all $b\in\bbR$ and if there exist two sequences $0<r_{1,n}<r_{2,n}<\infty$ such that

\begin{equation*} \begin{split} \begin{cases}
&\underset{u\in\mathscr{E}_{\leq n}, \vert u\vert=r_{2,n}}{\max}\;\Phi(u) \leq 0,   \\
&\underset{u\in\mathscr{E}_{>n}, \vert u\vert=r_{1,n}}{\inf}\;\Phi(u) \underset{n\rightarrow+\infty}{\longrightarrow}+\infty,
\end{cases} \end{split} \end{equation*}
then $\Phi$ has an unbounded sequence of critical values.
\end{thm}

\medskip

The parity condition on $f$ implies that $\Phi$ is even, hence invariant by the action of the multiplicative group $\{\pm 1\}$. The role played by the (no retraction)/(Brouwer fixed point) argument in the proof of Theorem \ref{linking} is played in that setting by the Borsuk-Ulam fixed point theorem. See e.g. Section 3.1 and Section 3.2 of \cite{Willem}.

\bigskip

\subsection{The Palais-Smale condition}
\label{SubsectionPalaisSmale}

We will work from now on with a nonlinearity $f\in C^1(\mcS\times\bbR,\bbR)$ that satisfies the following conditions, referred to in the text as \textbf{\textsf{Assumption (A)}}. Recall from \eqref{EqDefnF} the definition of $F$.

\medskip

\begin{itemize}
   \item[{\color{gray} $\bullet$}] {\it There is an exponent $\ell>2$ such that one has
$$
\big|f(x,z)\big|\lesssim 1+|z|^{\ell-1}, \qquad \big|\partial_zf(x,z)\big|\lesssim 1+|z|^{\ell-2},
$$
and $f(x,z)=o(z)$, as $z$ goes to $0$, uniformly in $x\in\mcS$.}   \hspace{0.15cm}
   
   \item[{\color{gray} $\bullet$}] {\it One has $F\geq 0$ and there exist $k>0$ and $\gamma>2$ such that for all $x\in\mcS$ one has 
   \begin{equation} \label{EqDifferentialConditionF}
   \gamma F(x,z)\leq zf(x,z),
   \end{equation} 
   on the set $\{|z|\geq k\}$.}
\end{itemize}
As an example, any focusing polynomial nonlinearity $f(x,z)=z^{2j+1}$ for an integer $j\ge 1$ satisfies Assumption (A).

\medskip

\begin{prop} \label{potential_PS}
The map $\Phi$ satisfies Palais-Smale condition $(\textsf{\emph{PS}})_b$ for all $b\in\bbR$.
\end{prop}

\medskip

\begin{Dem}
As a preliminary remark note that the differential condition \eqref{EqDifferentialConditionF} on the set $\{\vert z\vert>k\}$ gives the existence of positive constants $c_1,c_2$ such that one has the global lower bound 
\begin{equation} \label{EqInequalityF}
F(x,z) \geq c_1\vert z\vert^\gamma - c_2
\end{equation}
on all of $\mcS\times\bbR$. Recall from \eqref{EqDefnIndexm} the definition of the index $m$. Let now $(u_n)$ be a sequence of elements of $\mathscr{E}$ such that $\sup_n\Phi(u_n)=:M<+\infty$ and $\Phi'(u_n)$ tends to $0$. Write
$$
u_n =: y_n + y'_n \in \mathscr{E}_{\leq m} \oplus \mathscr{E}_{>m}.
$$
We will choose below a constant $\beta\in(\frac{1}{\gamma},\frac{1}{2})$. Independently of this constant, one has for $n$ large enough, say $n\geq n_0$, the inequality $\vert \Phi'(u_n)(v)\vert \leq \Vert v\Vert_\mathscr{E}$, for all $v\in\mathscr{E}$. One thus has for such indices
\begin{equation} \label{EqBound} \begin{split}
M+\|u_n\|_\mathscr{E} &\geq \Phi(u_n)-\beta \Phi'(u_n)(u_n)   \\
&= \left(\frac{1}{2}-\beta\right)\left(\|u_n\|_\mathscr{E} ^2 + \int_{\mcS}au_n^2\right)-\int_{\mcS}\big(F(\cdot,u_n)-\beta f(\cdot,u_n)u_n\big)   \\
&\geq \left(\frac{1}{2}-\beta\right)\left(\|u_n\|_\mathscr{E} ^2 + \int_{\mcS}au_n^2\right)+(\gamma\beta-1)\big(c_1\|u_n\|_{L^\gamma}^\gamma-c_2\big),
\end{split} \end{equation}
from Assumption (A) and \eqref{EqInequalityF}. Since the decomposition $u_n=y_n+y'_n$ is orthogonal in $L^2(\mcS)$ and the space $\mathscr{E}_{>m}$ is stable for the map $(-H_c+a)$ we can use the definition of $\mu_0$ and $\delta$ in Lemma \ref{LemDefnSmallDelta}, to get
\begin{align*}
\|u_n\|_\mathscr{E}^2 + \int_{\mcS} au_n^2&= \|y_n\|_\mathscr{E}^2 + \|y'_n\|_\mathscr{E}^2 + 2\langle y_n,y'_n\rangle_\mathscr{E} + \int_{\mcS}a\big(y_n^2+(y'_n)^2+2y_ny'_n\big)   \\
&=\left(\|y_n\|_\mathscr{E}^2 + \int_{\mcS} ay_n^2\right) + \left(\|y'_n\|_\mathscr{E}^2 + \int_{\mcS} a(y'_n)^2\right) + 2\underbrace{\left(\langle y_n, y'_n\rangle_\mathscr{E}+\int_{\mcS}ay_ny'_n\right)}_{=0}   \\
&\geq\mu_0\|y_n\|_{L^2}^2+\delta\|y'_n\|_\mathscr{E}^2.
\end{align*}
For the $L^\gamma$ norm in \eqref{EqBound} we remark that since $\mcS$ is compact and $\gamma>2$ the space $L^\gamma(\mcS)$ is a subspace of $L^2(\mcS)$ with
$$
\|u_n\|_{L^\gamma}^\gamma \apprge \|u_n\|_{L^2}^\gamma \apprge \left(\|y_n\|_{L^2}^2+\|y'_n\|_{L^2}^2\right)^{\gamma/2} \apprge \|y_n\|_{L^2}^\gamma.
$$
We thus have for $n$ large enough the inequality
$$
C+\|u_n\|_\mathscr{E} \geq \left(\frac{1}{2}-\beta\right) \Big(\mu_0\|y_n\|_{L^2}^2 + \delta\|y'_n\|_\mathscr{E}^2\Big) + c_1(\gamma\beta-1)\|y_n\|_{L^2}^\gamma
$$
for a positive constant $C$. Using the equivalence of the norms on the finite dimensional space $\mathscr{E}_{\leq m}$ where $y_n$ lives and choosing $\beta<1/2$ close enough to $1/2$ to have $\gamma\beta>1$ and the constant in front of $\|y_n\|_{L^2}^\gamma$ in 
$$
C+\|u_n\|_\mathscr{E} \geq \big(1/2-\beta\big)\delta\|y'_n\|_\mathscr{E}^2 + \Big(c_1(\gamma\beta-1)+ c_2\big(1/2-\beta\big)\mu_0\Big)\|y_n\|_{L^2}^\gamma
$$
positive, this implies that the sequence $u_n$ is bounded in $\mathscr{E}$. Indeed, assume for instance $\|y_n\|_{\mathscr{E}}$ is not bounded, then, as $\gamma>2$ the previous inequality rewrites

$$
1+\|y'_n\|_\mathscr{E}+\|y_n\|_{\mathscr{E}} \gtrsim \|y'_n\|_\mathscr{E}^2 + \|y_n\|_{\mathscr{E}}^2 \gtrsim \big(\|y'_n\|_{\mathscr{E}}+\|y_n\|_{\mathscr{E}}\big)^2
$$
proving that $\|y'_n\|_{\mathscr{E}}+\|y_n\|_{\mathscr{E}}$ is bounded, which contradicts the fact that $\|y_n\|_{\mathscr{E}}$ is not bounded. Similar argument ensures that $\|y'_n\|_{\mathscr{E}}$ is bounded as well, and thus so is $\|u_n\|_{\mathscr{E}}$.
  
\ssk

There is thus a subsequence $(u_{n'})$ that converges weakly to an element $u\in\mathscr{E}$ and in $L^{p/2}(\mcS)$ and $L^{2p/(p-1)}(\mcS)$ to $u$ as well. We obtain the convergence of $u_{n'}$ to $u$ in $\mathscr{E}$ from the identity
$$
\big(\Phi'(u_{n'})-\Phi'(u)\big)(u_{n'}-u) = \|u_{n'}-u\|_\mathscr{E}^2 - \int_{\mcS} \Big((f(\cdot,u_{n'}) - f(\cdot,u))(u_{n'}-u)-a(u_{n'}-u)^2\Big)
$$
and the fact that 
\begin{itemize}
   \item[--] the quantity $\big(\Phi'(u_{n'})-\Phi'(u)\big)(u_{n'}-u)$ is converging to $0$ since $u_{n'}$ is converging weakly to $u$ in $\mathscr{E}$ and $\Phi'(u_{n'})$ is converging to $0$,
   \item[--] the two quantities $\int_{\mcS}\big(f(\cdot,u_{n'})-f(\cdot,u)\big)(u_{n'}-u)$ and $\int_{\mcS}a(u_{n'}-u)^2$ are converging to $0$ from H\"older inequality and the $L^{p/2}(\mcS)$, respectively $L^{2p/(p-1)}(\mcS)$, convergence of $u_{n'}$ to $u$.
\end{itemize}
This concludes the proof that $\Phi$ satisfies the Palais-Smale condition $(\textsf{PS})_b$ for all $b\in\bbR$.
\end{Dem}

\bigskip

\subsection{Existence and multiplicity results}
\label{SubsectionExistenceMultiplicity}

We can now state and prove our main existence and multiplicity results for the semilinear equation
\begin{equation} \label{EqModelEqPrime}
-Hu+au=f(\cdot,u).
\end{equation}
Note that unlike in the fixed point approach of Proposition \ref{PropSchauder} no small size or a good integrability assumption on $a$ is needed in the next statement.

\medskip

\begin{thm} \label{link_PDE} 
If $f$ satisfies assumption \textbf{\textsf{(A)}}, then for any $a\in L^p(\mcS)$ with $p>1$, the equation \eqref{EqModelEqPrime} has a non-trivial weak solution in $\mathscr{E}$.
\end{thm}

\medskip

\begin{Dem}
As trading $a$ for $a-c$ does not change its integrability properties we consider the equation
$$
-H_cu+au=f(\cdot,u).
$$
Proposition \ref{potential_PS} shows that the map $\Phi$ satisfies the Palais-Smale condition $(\textsf{PS})_b$ for all $b\in\bbR$. We now check the condition \eqref{EqRabinowitzCondition} of Rabinowitz' linking theorem, Theorem \ref{linking}, with ${\sf y} = r_1\frac{e_{m+1}}{\Vert e_{m+1}\Vert_\mathscr{E}}$, for an appropriate choice of constants $0<r_1<r_2<\infty$. We use the notations of Theorem \ref{linking}.

\ssk

We have from the large and small $z$ behaviour of $f(x,z)$ stated in Assumption \textbf{\textsf{(A)}} the existence for any $\theta>0$ of a positive constant $c_\theta$ such that $|F(x,z)|\leq\theta |z|^2 + c_\theta|z|^\ell$, for all $(x,z)\in\mcS\times\bbR$. This gives in particular, for any $u\in\mathscr{E}_{>m}$, the lower bound
\begin{equation*} \begin{split}
\Phi(u) &\geq \frac{\delta}{2}\|u\|_\mathscr{E}^2 - \theta\|u\|_{L^2}^2 - c_\theta\|u\|_{L^\ell}^\ell   \\
 	   &\geq \Big(\frac{\delta}{2}-\theta\Big)\|u\|_\mathscr{E}^2 - c'_\theta\|u\|_{\mathscr{E}}^\ell,
\end{split} \end{equation*}
with $\delta$ as in Lemma~\ref{LemDefnSmallDelta}, and for another positive constant $c'_\theta$, from the embedding of $\mathscr{E}$ in $L^\ell(\mcS)$ when $\ell<2/(1-\alpha)$. As $\ell>2$ this inequality guarantees that for $0<\theta<\delta/2$ and $r_1$ small enough
$$
\underset{S_{r_1}\cap \mathscr{E}_{>m}}{\inf}\;\Phi > 0,
$$
with $S_{r_1}$ the sphere of $\mathscr{E}$ of radius $r_1$. We check in the sequel of the proof that one can find $r_2>r_1$ finite such that 
$$
\underset{\partial \mcB_{r_2}}{\sup}\;\Phi \leq 0.
$$
For $u\in\mathscr{E}_{\leq m}$ one has from the fact that $F$ is non-negative and $\mu_m$ non-positive
\begin{equation*} \begin{split}
\Phi(u)&=\frac{1}{2}\left(\|u\|_\mathscr{E}^2+\int_{\mcS}au^2\right) - \int_{\mcS}F(\cdot,u)   \\
&\leq\int_{\mcS}\left(\frac{\mu_m}{2}u^2-F(\cdot,u)\right) \leq 0.
\end{split} \end{equation*}
For any $r_2>0$ and $u=y+t{\sf y}\in\mcB_{r_2}$ we have from the global lower bound \eqref{EqInequalityF} on $F$, and the equivalence of norms on the finite dimensional space $\mathscr{E}_{\leq m}\oplus\bbR{\sf y}$, the estimate
\begin{equation*} \begin{split}
\Phi(u) &\leq\frac{1}{2}\|u\|_\mathscr{E}^2+\frac{1}{2}\|a\|_{L^p}\|u\|_{L^{2p/(p-1)}}^2-c_1\int_{\mcS}|u|^\gamma + c_2   \\
	   &\lesssim \|u\|_\mathscr{E}^2+1 - c'_1\|u\|_\mathscr{E}^\gamma,
\end{split} \end{equation*}
for a positive constant $c'_1$. It follows that $\Phi(u)\leq 0$ if $\Vert u\Vert_\mathscr{E}$ is large enough, since $\gamma>2$. The radius $r_2$ is chosen accordingly.
\end{Dem}

\medskip

The above proof males it clear that Theorem \ref{link_PDE} holds under the slightly weaker assumption that $F(\cdot,z)$ is only bounded below by $\mu_m z^2$.

\medskip

\begin{thm} \label{fountain_PDE} 
Assume in addition to assumption \textbf{\textsf{(A)}} that $f$ is odd with respect to its second argument. Then for any $a\in L^p(\mcS)$ for some $p>1$, there exists a sequence $(u_n)\subset\mathscr{E}$ of weak solutions of the equation
$$
Hu = au + f(\cdot,u)
$$
such that $\Phi(u_n)$ goes to $+\infty$ as $n$ goes to $+\infty$.
\end{thm}

\medskip

\begin{Dem}
We check that the conditions of Bartsch's fountain theorem (Theorem~\ref{ThmBartsch}) are met. Given $n\geq m$ and $u\in \mathscr{E}_{>n}$, for $\theta<\frac{\delta}{2}$, as in the proof of Theorem \ref{link_PDE}, we have
\begin{equation*} \begin{split}
\Phi(u)&=\frac{1}{2}\left(\|u\|_\mathscr{E}^2 + \int_{\mcS}au^2\right)-\int_{\mcS}F(\cdot,u)   \\
&\geq\frac{\delta}{2}\|u\|_\mathscr{E}^2 - \theta\|u\|_\mathscr{E}^2 - c_\theta\|u\|_{L^\ell}^\ell \geq\delta'\|u\|_\mathscr{E}^2 - c_\theta\beta_n^\ell\|u\|_\mathscr{E}^\ell   
\end{split} \end{equation*}
for a positive constant $\delta'$ and $\beta_n\defeq\underset{u\in \mathscr{E}_{>n}}{\sup}\frac{\|u\|_{L^\ell}}{\|u\|_\mathscr{E}}$. Set 
$$
r_{1,n}^{\ell-2} \defeq \frac{\delta'}{2c_\theta\beta_n^\ell}
$$ 
and take any $u\in\mathscr{E}_{>n}$ with $\|u\|_\mathscr{E} = r_{1,n}$. Then we have
$$
\Phi(u)\geq r_{1,n}^2\left(\delta'-c_\theta\beta_n^\ell r_{1,n}^{\ell-2}\right) = \frac{\delta'}{2}r_{1,n}^2
$$
In turns out that $r_{1,n}$ diverges to $+\infty$. To see this, note that the $\beta_n$ are non-increasing so they have a limit $\beta\geq 0$. Pick for each $n\geq m$ a point $u_n\in \mathscr{E}_{>n}$ such that $\|u_n\|_\mathscr{E}=1$ and $\|u_n\|_{L^\ell}\geq \beta_n/2$. Up to extraction, the $u_n$ are converging weakly in $\mathscr{E}$ and in $L^\ell(\mcS)$ to a limit element $u\in\mathscr{E}$. But it follows from the definition of $\mathscr{E}_{>n}$ that the  $u_n$ are converging weakly to $0$, so $\beta=0$ and $r_{1,n}$ diverges to $+\infty$.

\ssk

To control the behaviour of $\Phi$ on $\mathscr{E}_{\leq n}$ we proceed as in the proof of Theorem \ref{link_PDE} and write for $u\in\mathscr{E}_{\leq n}$

\begin{equation*} \begin{split}
\Phi(u)&= \frac{1}{2}\left(\|u\|_\mathscr{E}^2+\int_{\mcS}au^2\right) - \int_{\mcS}F(\cdot,u)   \\
&\leq \frac{1}{2}\|u\|_\mathscr{E}^2+\frac{1}{2}\|a\|_{L^p}\|u\|_{L^{2p/(p-1)}}^2 - c_1\int_{\mcS}|u|^\gamma + c_2   \\
&\leq C_1(\|u\|_\mathscr{E}^2+1)-C_{2,n}\|u\|_\mathscr{E}^\gamma, 
\end{split} \end{equation*}
for some positive constant $C_1$ and some $n$-dependent constant $C_{2,n}$, using the equivalence of norms on the finite dimensional space $\mathscr{E}_{\leq n}$. The condition $\gamma>2$ thus guarantees that $\Phi$ takes non-positive values on the intersection with $\mathscr{E}_{\leq n}$ of the sphere of $\mathscr{E}$ of a well-chosen radius $r_{2,n}>r_{1,n}$.
\end{Dem}

\medskip

\begin{cor}
For any non-null even integer $\ell$ and any potential $a\in L^p(\mcS)$, with $p>1$, the semilinear problem
$$
-Hu+au = u|u|^\ell
$$
has infinitely many weak solutions.
\end{cor}

\bigskip

\section{A non-variational singular stochastic PDE}
\label{SectionNonVariational}

We consider in this section the case of the two dimensional torus $\mcS=\bbT^2$. Denote by $\star$ the convolution operation in $\bbT^2$. We consider in this section the singular Choquard-Pekar equation
\begin{equation} \label{EqChoquardPekar}
(-H_c +a)u = (w\star |u|^p)|u|^{q-2}u,
\end{equation}
for appropriate parameters $w$, $p$, and $q$, which can be seen as a generalization of the (stationary) Hartree equation on $\bbT^2$ (for a survey on these equations in the deterministic case and the corresponding parameters, see \cite{Ackermann}). While the latter can be treated with variational methods, \eqref{EqChoquardPekar} cannot be written as the Euler-Lagrange equation of a functional on $\scrE$ as soon as $p\neq q$, the case of interest here. We use Ghoussoub's machinery of self-dual functionals to tackle that equation. We recall what we need from this setting in the restricted functional setting of the space $\mathscr{E}$ -- this will be sufficient for us. See Ghoussoub's book \cite{Ghoussoub} for the whole story. It will clarify things here to make a difference between the Hilbert space $\scrE$ and its topological dual $\scrE'$ without identifying the later to the former.

\medskip

Given a convex and lower semi continuous functional $\varphi : \scrE \to\IR$ its Fenchel transform $\varphi' : \scrE'\rightarrow \IR$ is defined by
$$
\varphi'(p)\defeq\sup_{u\in\scrE} \big(p(u) - \varphi(u)\big),
$$
and its subdifferential at a point $u\in\scrE$ is the subset of $\scrE'$ defined by
$$
\partial\varphi(u) \defeq \Big\{p\in \scrE'\,;\, \forall h\in \scrE, \varphi(u+h)\geq \varphi(u) + p(h)\Big\}.
$$
One thus has
$$
\varphi(u)+\varphi'(p)\geq p(u)
$$
for all $u\in\scrE, p\in\scrE'$, with equality if and only if $p\in\partial\varphi(u)$. An operator $\Lambda : \scrE\to\scrE'$ is said to be {\it regular} if it is weak-to-weak continuous on its domain and $u\mapsto (\Lambda u)(u)$ is weakly lower semi-continuous. Recall also that a non-negative function $\Phi : \scrE \to[0,+\infty)$ is said to be {\it self-dual} if there exists a real-valued function $M$ on $\scrE\times\scrE$ such that 
\begin{equation} \label{EqPhiFromM}
\Phi(u) = \sup M(u,\cdot)
\end{equation}
for all $u\in\scrE$, where all the functions $M(v,\cdot)$ are proper and concave, all the functions $M(\cdot,v)$ are weakly semicontinuous and $M$ is non-positive on the diagonal. A large class of self-dual functions is provided in the following statement. It is a direct consequence of Theorem 12.3 in Ghoussoub's book \cite{Ghoussoub} -- itself a direct consequence of Ky Fan's min-max principle.

\medskip

\begin{thm} \label{ThmGhoussoub}
Let $\varphi:\scrE\to \bbR$ be a lower semicontinuous convex function that is bounded below. Let $f\in\scrE'$ and $\Lambda : \scrE\rightarrow\scrE'$ be a regular (possibly nonlinear) operator. Then the function
$$
M(u,v) := (\Lambda u)(u-v) + \varphi(u) - \varphi(v)
$$
on $\scrE\times\scrE$ defines via \eqref{EqPhiFromM} a non-negative self-dual functional 
$$
\Phi(u) = \varphi(u) - \varphi'(\Lambda u) + (\Lambda u)(u).
$$
If further $\big\Vert \varphi(u) + (\Lambda u)(u)\big\Vert_{\scrE}$ tends to $+\infty$ as $\|u\|_{\scrE}$ tends to $+\infty$ then the function $\Phi$ attains its minimum $0$ at some point $\overline{u}\in\scrE$ where
$$
-\Lambda \overline{u}\in\partial\varphi(\overline{u}).
$$
\end{thm}

\medskip

We will use Theorem \ref{ThmGhoussoub} to prove the next statement, with $\varphi$ the $C^1$ function
$$
\varphi(u) := \frac{1}{2}\|u\|_{\scrE}^2+\frac{1}{2}\int_{\bbT^2} au^2
$$
on $\scrE$. Its subdifferential being a singleton the conclusion of Theorem \ref{ThmGhoussoub} will thus come under the form that $\overline{u}$ is a weak solution of the equation
$$
\partial\varphi(u) + \Lambda u = 0, 
$$
that is
$$
(-H_c+a)\overline{u} + \Lambda \overline{u} = 0.
$$
An ad hoc choice of function $\Lambda$ will identify this equation with the Choquard-Pekar equation \eqref{EqChoquardPekar}.

\medskip

\begin{thm} \label{ThmChoquardPekar}
Pick exponents $p\in[1,+\infty)$, $q\in(1,+\infty)$ and let the potential $a$ be bounded and positive. Assume that the interaction kernel $w\in L^1(\bbT^2)$ is non-positive. Then the singular Choquard-Pekar equation \eqref{EqChoquardPekar} has a weak solution $\overline{u}\in\scrE$.
\end{thm}

The non-positivity assumption on the interaction kernel $w$ may seem rather ad hoc, as the Choquard-Pekar equation originally arose in the physics litterature for a non-\textit{negative} kernel. However it also corresponds to the variational case $p=q$, and as we pointed out above, we rather see the equation \eqref{EqChoquardPekar} as a toy-model to extend the self-dual machinery to the case of a singular stochastic PDE. From this point of view our Theorem \ref{ThmChoquardPekar} is really the analogue of \cite[Theorem 12.5 (A)]{Ghoussoub}.
\medskip

\begin{Dem}
The boundedness and positivity assumption on the potential $a$ guarantees that the function $\varphi$ on $\scrE$ is well-defined, convex, non-negative and lower semicontinuous. For $u\in\scrE$ set
$$
\Lambda u \defeq -(w\star |u|^p)|u|^{q-2}u.
$$
One has for all $u,v\in\scrE$
\begin{equation*} \begin{split}
\big|(\Lambda u)(v)\big| &\leq \int_{\IT^2}\big|w\star|u|^p\big|\, |u|^{q-1}|v|   \\
					       &\leq \|w\|_{L^1} \, \||u|^p\|_{L^2} \, \big\||u|^{q-1}v\big\|_{L^2}   \\
					       &\leq \|w\|_{L^1} \, \|u\|_{L^{2p}}^p \, \|u\|_{L^{2q}}^{q-1} \, \|v\|_{L^{2q}},
\end{split} \end{equation*}
where we used Hölder inequality to bound $\big\||u|^{q-1}v\big\|_{L^2}$. Recall the compact embedding of $\scrE$ into $L^r(\bbT^2)$ for any $1<r<\frac{2}{1-\alpha}$ mentioned in section 2.1 (where $\alpha-2$ stands for the regularity of the white noise $\xi$). Using it for $r\in\{2p,2q\}$ and an appropriate choice of $\alpha$ close enough to 1, then yields the bound
$$
\big|(\Lambda u)(v)\big| \lesssim \|w\|_{L^1}\|u\|_\scrE^{p+q-1}\|v\|_\scrE
$$
that shows that $\Lambda$ is a well-defined map from $\scrE$ into $\scrE'$.

We check the weak-to-weak continuity of $\Lambda$. Let $(u_n)$ converge weakly to $u$ in $\scrE$. Let $v\in\scrE$. We have
\begin{equation*} \begin{split}
\big| (\Lambda u_n)(v) - (\Lambda u)(v) \big| &\leq \left|\int_{\bbT^2} \big(w\star |u_n|^p\big)|u_n|^{q-2}u_nv - \int_{\bbT^2} \big(w\star |u|^p\big) |u_n|^{q-2}u_nv\right|   \\
&\qquad+\left|\int_{\bbT^2} \big(w\star |u|^p\big) |u_n|^{q-2}u_nv - \int_{\bbT^2} \big(w\star |u|^p\big) |u|^{q-2}uv\right|   \\
&\leq\left|\int_{\bbT^2} \Big(w\star\big(|u_n|^p-|u|^p\big)\Big)|u_n|^{q-2}u_nv\right|   \\
&\qquad+ \left|\int_{\bbT^2} \big(w\star |u|^p\big)\Big(|u_n|^{q-2}u_n-|u|^{q-2}u\Big)v \right|
\end{split} \end{equation*}
One can use once again the compact embedding of $\scrE$ into $L^r(\bbT^2)$ for some appropriate choice of $\alpha<1$ and the weak convergence of $u_n$ to 
$u$ to get the convergence of $|u_n|^p$ to $|u|^p$ in $L^2(\bbT^2)$ and the convergence of $|u_n|^{q-2}u_n$ to $|u|^{q-2}u$ in $L^{2q/(q-1)}(\bbT^2)$. We therefore have

\begin{align*}
\big| (\Lambda u_n)(v) - (\Lambda u)(v) \big| \lesssim &\|v\|_{L^{2q}}\big\||u_n|^p-|u|^p\big\|_{L^2}\|u_n\|_{L^{2q}}^{q-1}   \\
&\qquad+ \|v\|_{L^{2q}}\||u|^p\|_{L^2} \big\||u_n|^{q-2}u_n-|u|^{q-2}u\big\|_{L^{2q/(q-1)}}
\end{align*}
with an implicit multiplicative constant depending on $\|w\|_{L^1}$. The upper bound vanishes as $n$ goes to $\infty$, which shows the weak-to-weak continuity. It follows in particular from the previous estimates that 
$$
\big| (\Lambda u_n)(u_n) - (\Lambda u)(u) \big| \leq \|w\|_{L^1} \big\||u_n|^p-|u|^p\big\|_{L^2}\||u_n|^q\|_{L^2} + \|w\|_{L^1}\||u|^p\|_{L^2} \big\||u_n|^q-|u|^q\big\|_{L^2}
$$
goes to $0$ as $n$ goes to $\infty$. All this proves that the function $\Lambda$ is regular. Remark at last that since the interaction kernel $w$ is non-positive the function $\varphi(u) + (\Lambda u) u$ is coercive. We are thus in the setting of Theorem \ref{ThmGhoussoub}, from which we get our conclusion.
\end{Dem}

\medskip

We note from the fact that $\Lambda$ takes its values in $\scrE'$ that we could try and use a fixed point strategy to get a solution of equation \eqref{EqChoquardPekar}, as in Proposition \ref{PropSchauder} or the comment following it. Assuming $a\in L^2(\bbT^2)$ and using Cameron-Martin theorem gives the existence of a random constant $c$ such that equation \eqref{EqChoquardPekar} has almost surely a solution. In any case this would require that we assume either that $a$ is small enough in $L^\infty(\bbT^2)$ or sufficiently integrable, and that $w$ is small enough in $L^1(\bbT^2)$. The use of Theorem \ref{ThmGhoussoub} bypasses this kind of constraints. It is straightforward to adapt the proof of Theorem \ref{ThmChoquardPekar} to the more general case of the equation 
$$
(-H_c +a)u = \big(w\star f(u)\big)g(u),
$$
for nonlinearities $f : \bbR\rightarrow\bbR$ and $g : \bbR\rightarrow\bbR$ that are uniformly continuous and such that $\big\vert f(z)\big\vert \lesssim 1+\vert z\vert^p$ and $\big\vert g(z)\big\vert \lesssim 1+\vert z\vert^{q-1}$. So the existence result of Theorem \ref{ThmChoquardPekar}  holds in that setting.

\bigskip
\bigskip

\noindent \textcolor{gray}{$\bullet$} {\sf I. Bailleul} -- Univ. Rennes, CNRS, IRMAR - UMR 6625, F-35000 Rennes, France.   \\
\noindent {\it E-mail}: ismael.bailleul@univ-rennes1.fr   

\medskip

\noindent \textcolor{gray}{$\bullet$} {\sf H. Eulry} --  Univ. Rennes, CNRS, IRMAR - UMR 6625, F-35000 Rennes, France.   \\
\noindent {\it E-mail}: hugo.eulry@ens-rennes.fr 

\medskip

\noindent \textcolor{gray}{$\bullet$} {\sf T. Robert} --  IECL -- Facult\'e des sciences et Technologies, 54506 Vandoeuvre-l\`es-Nancy, France.   \\
\noindent {\it E-mail}: tristan.robert@univ-lorraine.fr 

\end{document}